\magnification=1200
\vsize=20.6cm \voffset=-2mm \hsize=14.1cm \hoffset=-4mm
\parindent=5mm
\parskip=4pt plus 2pt minus 1pt

\font\hugebf=cmbx10 at 14.4pt
\font\bigbf=cmbx10 at 12pt
\font\smallrm=cmr10 at 8pt
\font\smallbf=cmbx10 at 8pt

\font\tenCal=eusm10
\font\sevenCal=eusm7
\font\fiveCal=eusm5
\newfam\Calfam
  \textfont\Calfam=\tenCal
  \scriptfont\Calfam=\sevenCal
  \scriptscriptfont\Calfam=\fiveCal
\def\Cal{\fam\Calfam\tenCal}

\font\tenBbb=msbm10
\font\sevenBbb=msbm7
\font\fiveBbb=msbm5
\newfam\Bbbfam
  \textfont\Bbbfam=\tenBbb
  \scriptfont\Bbbfam=\sevenBbb
  \scriptscriptfont\Bbbfam=\fiveBbb
\def\Bbb{\fam\Bbbfam\tenBbb}

% Usual sets of numbers  
\def\bC{{\Bbb C}}
\def\bN{{\Bbb N}}
\def\bP{{\Bbb P}}

\def\bR{{\Bbb R}}

\def\cO{{\Cal O}}
\def\cS{{\Cal S}}

\def\rank{\mathop{\rm rank}}
\def\Id{{\rm Id}}
\def\tr{\mathop{\rm tr}\nolimits}
\def\Hom{\mathop{\rm Hom}\nolimits}
\def\Herm{\mathop{\rm Herm}\nolimits}
\def\MAVol{\mathop{\rm MAVol}\nolimits}
\def\det{\mathop{\rm det}\nolimits}
\def\bu{\hbox{$\scriptstyle\bullet$}}

\let\ol=\overline
\def\square{{\hfill \hbox{
\vrule height 1.453ex  width 0.093ex  depth 0ex
\vrule height 1.5ex  width 1.3ex  depth -1.407ex\kern-0.1ex
\vrule height 1.453ex  width 0.093ex  depth 0ex\kern-1.35ex
\vrule height 0.093ex  width 1.3ex  depth 0ex}}}
\def\qed{$\square$}

\def\Langle{\langle\!\langle}
\def\Rangle{\rangle\!\rangle}

\def\Bibitem#1&#2&#3&#4&%
{\hangindent=1.6cm\hangafter=1
\noindent\rlap{\hbox{\bf #1}}\kern1.6cm{\rm #2}{\it #3}{\rm #4.}} 

\strut\vskip1cm
\centerline{\hugebf Hermitian-Yang-Mills approach to the conjecture of}\vskip6pt
\centerline{\hugebf Griffiths on the positivity of ample vector bundles}
\vskip5mm
\centerline{\bf Jean-Pierre Demailly$\,{}^{\bf*}$ \footnote{}{\sevenrm
\baselineskip=9pt\strut\kern-14.5pt* This work is supported by the
European Research Council project ``Algebraic and K\"ahler Geometry''
\hfil\break
\hbox{(ERC-ALKAGE, grant No. 670846 from September 2015)}\vskip-9pt}}

\centerline{Universit\'e Grenoble Alpes, Institut Fourier}
\centerline{March 17, 2020, revised February 4, 2021}
\vskip4mm
\line{\hfill\it Dedicated to Professor Vyacheslav Shokurov on the occasion
of his 70${}^{\,th}$ birthday}
\vskip4mm

{\noindent\smallrm\baselineskip=9.6pt\par\noindent
{\smallbf Abstract.} Given a vector bundle of arbitrary rank with ample
determinant line bundle on a projective manifold, we propose a new
elliptic system of differential equations of Hermitian-Yang-Mills type
for the curvature tensor. The system is designed so that
solutions provide Hermitian metrics with positive curvature in the sense
of Griffiths~-- and even in the dual Nakano sense. As a consequence,
if an existence result could be obtained for every ample vector bundle,
the Griffiths  conjecture on the equivalence between ampleness and
positivity of vector bundles would be settled.
\medskip

\noindent
{\smallbf Keywords.} Ample vector bundle, Hermitian metric, Griffths positivity,
Nakano positivity, Hermitian-Yang-Mills equation, Monge-Amp\`ere equation,
elliptic operator.
\medskip

\noindent
{\smallbf MSC Classification 2020.} 32J25, 53C07
\vskip4mm}

\noindent
{\bigbf 1. Introduction}\medskip

Let $X$ be a projective $n$-dimensional manifold. A conjecture due to
Griffiths [Gri69] stipulates that a holomorphic vector bundle $E\to X$
is ample in the sense of Hartshorne, meaning that the associated
line bundle $\cO_{\bP(E)}(1)$ is ample, if and only if $E$ possesses a
Hermitian metric $h$ such that the Chern curvature tensor
$\Theta_{E,h}=i\nabla_{E,h}^2$ is Griffiths positive. In other words,
if we let $\rank E=r$ and
$$
\Theta_{E,h}=i\sum_{1\le j,k\le n,\,1\le\lambda,\mu\le r}c_{jk\lambda\mu}dz_j\wedge
d\ol z_k\otimes e_\lambda^*\otimes e_\mu\leqno(1.1)
$$
in terms of holomorphic coordinates $(z_1,\ldots,z_n)$ on $X$ and of
an orthonormal frame $(e_\lambda)_{1\le\lambda\le r}$ of $E$,
the associated quadratic form
$$
\widetilde\Theta_{E,h}(\xi\otimes v):=
\langle\Theta_{E,h}(\xi,\ol\xi)\cdot v,v\rangle_h=
\sum_{1\le j,k\le n,\,1\le\lambda,\mu\le r}
c_{jk\lambda\mu}\xi_j\ol\xi_k v_\lambda\ol v_\mu
\leqno(1.2)
$$
should take positive values on non zero tensors $\xi\otimes v\in T_X\otimes E$.
A stronger concept is Nakano positivity (cf.\ [Nak55]), asserting that
$$
\widetilde\Theta_{E,h}(\tau):=\sum_{1\le j,k\le n,\,1\le\lambda,\mu\le r}
c_{jk\lambda\mu}\tau_{j\lambda}\ol\tau_{k\mu}>0
\leqno(1.3)
$$
for all non zero tensors $\tau=\sum_{j,\lambda}
\tau_{j\lambda}{\partial\over\partial z_j}\otimes e_\lambda\in T_X\otimes E$.
It is in fact interesting to consider also the curvature tensor of the dual
bundle $E^*$, which happens to be given by the opposite of the transpose
of $\Theta_{E,h}$, i.e.\
$$
\Theta_{E^*,h^*}=-{}^T\Theta_{E,h}=-
\sum_{1\le j,k\le n,\,1\le\lambda,\mu\le r}c_{jk\mu\lambda}dz_j\wedge
d\ol z_k\otimes (e^*_\lambda)^*\otimes e^*_\mu.
\leqno(1.4)
$$
This leads to the concept of dual Nakano positivity, stipulating that
$$
-\widetilde\Theta_{E^*,h^*}(\tau)=
\sum_{1\le j,k\le n,\,1\le\lambda,\mu\le r}
c_{jk\mu\lambda}\tau_{j\lambda}\ol\tau_{k\mu}>0
\leqno(1.5)
$$
for all non zero tensors $\tau=\sum_{j,\lambda}
\tau_{j\lambda}{\partial\over\partial z_j}\otimes e^*_\lambda\in T_X\otimes E^*$.
On the other hand, Griffiths positivity of $\Theta_{E,h}$ is equivalent to 
Griffiths negativity of $\Theta_{E^*,h^*}$, and implies the positivity of the
induced metric on the tautological line bundle $\cO_{\bP(E)(1)}$.
By the Kodaira embedding theorem [Kod54], the positivity of $\cO_{\bP(E)(1)}$
is equivalent to its ampleness, hence we see immediately from the definitions
that
$$
\widetilde\Theta_{E,h}~\hbox{dual Nakano positive}~~\Rightarrow~~
\widetilde\Theta_{E,h}~\hbox{Griffiths positive}~~\Rightarrow~~
E~\hbox{ample}.\leqno(1.6)
$$
In this short note, we consider the following converse problem:\vskip8pt

\vbox{\noindent
{\bf 1.7. Basic question.} {\it Does it hold that}\vskip7pt
\line{\hfill $E~\hbox{ample}~\Rightarrow~
\widetilde\Theta_{E,h}~\hbox{dual Nakano positive}~?$\hfill}}
\medskip

\noindent A positive answer would clearly settle the Griffiths
conjecture, in an even stronger form. One should observe
that Nakano positivity implies Griffiths positivity, but is in
general a more restrictive condition. As a consequence, one cannot
expect ampleness to imply Nakano positivity. For instance, $T_{\bP^n}$
is easily shown to be ample (and Nakano semi-positive for the
Fubini-Study metric), but it is not Nakano positive, as the Nakano
vanishing theorem [Nak55] would then yield
$$H^{n-1,n-1}(\bP^n,\bC)=
H^{n-1}(\bP^n,\Omega^{n-1}_{\bP^n})=
H^{n-1}(\bP^n,K_{\bP^n}\otimes T_{\bP^n})=0.
\leqno(1.8)
$$
On the other hand, it does not seem that there are any examples of
ample vector bundles that are not dual Nakano positive, thus the above basic
question is still legitimate, even though it might look very optimistic. 
We should mention here that subtle
relations between ampleness, Griffiths and Nakano positivity
are known to hold -- for instance, B.~Berndtsson [Ber09] has proved that
the ampleness of $E$ implies the Nakano positivity of $S^mE\otimes\det E$ for
every~$m\in\bN$. See also [DeS79] for an earlier direct and elementary proof
of the much weaker result that the Griffiths positivity of $E$
implies the Nakano positivity of $E\otimes\det E$, and [MoT07]
for further results analogue to those of~[Ber09].

So far, the Griffiths conjecture is known to hold when $n=\dim X=1$ or
$r=\rank E=1$ (in which cases, Nakano and dual Nakano positivity
coincide with Griffiths positivity). Proofs can be found in [Ume73,
Theorem 2.6] and [CaF90]. In both cases, the proof is based on the
existence of Harder-Narasimhan filtrations and on the 
Narasimhan-Seshadri theorem [NaS65] for stable vector bundles --
the 1-dimensional case of the Donaldson-Uhlenbeck-Yau theorem [Don85], [UhY86].
It is tempting to investigate whether techniques of gauge theory could be
used to approach the Griffiths conjecture. In this direction, P.~Naumann
[Nau17] proposed a K\"ahler-Ricci flow method that starts
with a given Finsler metric of positive curvature, and converges to
a Hermitian metric. It is however unclear whether the flow introduced
in [Nau17] preserves positivity, so it might very well produce in the limit
a Hermitian metric that does not have positive curvature. Another related
suggestion is V.~Pingali's proposal made in [Pin20] to study the vector
bundle Monge-Amp\`ere equation $(\Theta_{E,h})^n=\eta\,\Id_E$, where
$\eta$ is a positive volume form on $X$. Solving such an equation requires
polystability in dimension $n=1$, and, in general, a positivity
property of $(E,h)$ that is even stronger than Nakano positivity (and
thus much stronger than ampleness).

In section~2, we describe a more flexible differential system based on
a combination of a huge determinantal equation and a trace free
Hermite-Einstein condition. It relies on the well known continuity
method, and is designed to enforce positivity of
the curvature, actually in the dual Nakano sense -- a condition that
could eventually still be equivalent to ampleness.
We show that it is possible to design a non linear differential system
that is elliptic and invertible, at least near the origin of time.
It would however remain 
to check whether  one can obtain long time existence of the solution
for the said equation or  one of its variants. Section~3 is devoted to 
the discussion of a related extremal problem, and a concept of
volume for vector bundles.\bigskip

\noindent
{\bigbf 2. Approach via a combination of Monge-Amp\`ere and\vskip0pt
Hermitian-Yang-Mills equations}\medskip

Let $E\to X$ be a holomorphic vector bundle equipped with a smooth Hermitian
metric~$h$. If the Chern curvature tensor $\Theta_{E,h}$ is dual
Nakano positive, then the ${1\over r}$-power of
the $(n\times r)$-dimensional determinant of the corres\-ponding
Hermitian quadratic form on $T_X\otimes E^*$ can be seen as a
positive $(n,n)$-form
$$
\det_{T_X\otimes E^*}({\,}^T\Theta_{E,h})^{1/r}
:=\det(c_{jk\mu\lambda})_{(j,\lambda),(k,\mu)}^{1/r}\,
i dz_1\wedge d\ol z_1\wedge\ldots\wedge i dz_n\wedge d\ol z_n.
\leqno(2.1)
$$
Moreover, this $(n,n)$-form does not depend on the choice of coordinates $(z_j)$
on~$X$, nor on the choice of the orthonormal frame $(e_\lambda)$ on $E$
(but the orthonormality of $(e_\lambda)$ is required). Conversely,
given a K\"ahler metric $\omega_0$ on $X$, the basic
idea is that assigning a ``matrix Monge-Amp\`ere equation''
$$
\det_{T_X\otimes E^*}({\,}^T\Theta_{E,h})^{1/r}=f\,\omega_0^n,
\leqno(2.2)
$$
where $f$ is a smooth positive function, may enforce the dual
Nakano positivity of $\Theta_{E,h}$
if that assignment is combined with a
continuity technique from an initial starting point where posi\-tivity
is known. For $r=1$, we have
${}^T\Theta_{E,h}=\Theta_{E,h}=-i\partial\ol\partial\log h$,
and equation (2.2) is a standard Monge-Amp\`ere equation. If $f$ is
given and independent of $h$, Yau's theorem [Yau78] guaran\-tees
the existence of a unique solution $\theta=\Theta_{E,h}>0$, provided
$E$ is an ample line bundle and $\int_Xf\,\omega_0^n=c_1(E)^n$. One then gets
a smoothly varying solution $\theta_t=\Theta_{E,h_t}>0$ when the right
hand side $f_t$ of (2.2) varies smoothly with respect to some parameter~$t$.

Now, assuming $E$ to be ample of rank $r>1$, equation (2.2) becomes
underdetermined, since the real rank of the space of hermitian matrices
$h$ on $E$ is equal to $r^2$, while (2.2) provides only one scalar equation.
If $E=\bigoplus_{1\le j\le r}E_j$ splits as a direct sum of ample line
bundles and we take a diagonal Hermitian structure
$h=\bigoplus h_j$ on $E$, the $nr\times nr$ determinant splits as
a product of blocks, and equation (2.2) reduces to
$$
\Bigg(\prod_{1\le j\le r}\Theta_{E_j,h_j}^n\Bigg)^{1/r}=f\,
\omega_0^n.
\leqno(2.2_s)
$$
This ``split equation'' can be solved for any $f=\prod f_j^{1/r}$ with
$\int_X f_j\,\omega_0^n=c_1(E_j)^n$, just by solving the
individual equations $\Theta_{E_j,h_j}^n=f_j\,\omega_0^n$,
$f_j>0$, but the decomposition need not be unique. In this case,
the H\"older inequality requires
$\int_X f\omega_0^n\le(\prod c_1(E_j)^n)^{1/r}$, and the equality
can be reached by taking all $f_j$'s to be proportional to~$f$.

In general, solutions might still exist, but the lack of uniqueness
prevents us from getting a priori bounds. In order to recover a well
determined system
of equations, one needs to introduce $(r^2-1)$ additional scalar equations,
or rather a matrix equation of real rank $(r^2-1)$. If $E$ is ample,
the determinant line bundle $\det E$ is also~ample. By the Kodaira
embedding theorem, we can find a smooth Hermitian metric $\eta_0$ on
$\det E$ so that $\omega_0:=\Theta_{\det E,\eta_0}>0$
is a K\"ahler metric on $X$. In case $E$ is $\omega_0$-stable
or \hbox{$\omega_0$-polystable}, we know by the Donaldson-Uhlenbeck-Yau
theorem that there exists a Hermitian metric $h$ on $E$ satisfying the
Hermite-Einstein condition
$$
\omega_0^{n-1}\wedge\Theta_{E,h}={1\over r}\,\omega_0^n\otimes \Id_E,
\leqno(2.3)
$$
since the slope of $E$ with respect to $\omega_0\in c_1(E)$ is equal
to ${1\over r}$.

In general, one cannot expect $E$ to be $\omega_0$-polystable,
but Uhlenbeck-Yau have shown that there always exist
smooth solutions to a certain ``cushioned'' Hermite-Einstein
equation. To~make things more precise, let $\Herm(E)$ be the space
of Hermitian (non necessarily positive) forms on~$E$, and
given a Hermitian metric $h>0$, let $\Herm_h(E,E)$ be the space of
$h$-Hermitian endomorphisms $u\in \Hom(E,E)$; we denote by
$$
\Herm(E)\to \Herm_h(E,E),\quad q\mapsto\widetilde q~~\hbox{such that}~~
q(v,w)=\langle v,w\rangle_q=\langle\kern1pt\widetilde q\kern1pt(v),w\rangle_h
\leqno(2.4)
$$
the natural isomorphism between Hermitian quadratic forms and Hermitian endomorphisms, which depends of course on $h$. We also let
$$
\Herm^\circ_h(E,E)=\big\{u\in\Herm_h(E,E)\,;\;\tr u=0\big\}
\leqno(2.5)
$$
be the subspace of ``trace free'' Hermitian endomorphisms.
In the sequel, we fix a reference
Hermitian metric $H_0$ on $E$ such that $\det H_0=\eta_0$, so~that
$\Theta_{\det E,\det H_0}=\omega_0>0$. By [UhY86, Theorem~3.1], for
every $\varepsilon>0$, there exists a smooth Hermitian metric $q_\varepsilon$
on $E$ such that
$$
\omega_0^{n-1}\wedge\Theta_{E,q_\varepsilon}=\omega_0^n\otimes
\bigg({1\over r}\Id_E-\varepsilon\,\log\widetilde q_\varepsilon\bigg),
\leqno(2.6)
$$
where $\widetilde q_\varepsilon$ is computed with respect to~$H_0$,
and $\log u$ denotes the logarithm of a positive Hermitian endomorphism~$u$.
The intuitive reason is that the term $\log\widetilde q_\varepsilon$ introduces
sufficient ``friction'' to avoid any explosion of approximating solutions
when using a standard continuity method (see sections 2,3 in [UhY86]). 
On the other hand, when $\varepsilon\to 0$, the
metrics $q_\varepsilon$  become ``more and more distorted'' and
yield asymptotically a splitting of $E$ in weakly holomorphic subbundles
corresponding to the Harder-Narasimhan filtration of $E$ with respect
to~$\omega_0$. If we write $\det q_\varepsilon=e^{-\varphi}\,\det H_0$ and
take the trace in (2.6), we find
$\omega_0^{n-1}\wedge(\omega_0+i\partial\ol\partial\varphi)=
\omega_0^n(1+\varepsilon\varphi)$, hence $\omega_0^{n-1}\wedge
i\partial\ol\partial\varphi-\varepsilon\varphi\,\omega_0^n=0$.
A~standard application of the maximum principle shows that $\varphi=0$,
thus (2.6) implies $\det q_\varepsilon=\det H_0$ and
$\log\widetilde q_\varepsilon\in\Herm^\circ_{H_0}(E,E)$. In general,
for an arbitrary Hermitian metric $h$, we let
$$
\Theta_{E,h}^\circ=\Theta_{E,h}-
{1\over r}\Theta_{\det E,\det h}\otimes\Id_E\in
C^\infty(X,\Lambda_\bR^{1,1}T^*_X\otimes\Herm_h^\circ(E,E))
\leqno(2.7)
$$
be the curvature tensor of $E\otimes(\det E)^{-1/r}$ with respect to the
trivial determinant metric $h^\circ:=h\otimes(\det h)^{-1/r}$. Equation (2.6) is
equivalent to prescribing
$\smash{\det q_\varepsilon}=\det H_0$ and
$$
\omega_0^{n-1}\wedge\Theta^\circ_{E,q_\varepsilon}=
-\varepsilon\,\omega_0^n\otimes\log\widetilde q_\varepsilon.
\leqno(2.8)
$$
This is a matrix equation of rank $(r^2-1)$
that involves only $q_\varepsilon^\circ$ and does not depend on
$\det q_\varepsilon$. Notice that we have here
$\log\widetilde q_\varepsilon\in\Herm_{H_0}^\circ(E,E)$, but also
$\log\widetilde q_\varepsilon\in\Herm_{q_\varepsilon}^\circ(E,E)$.

In this context, given $\alpha>0$ large enough,
it~seems natural to search for a time dependent family of
metrics $h_t(z)$ on the fibers $E_z$ of $E$, $t\in [0,1]$, satisfying a 
generalized Monge-Amp\`ere equation
$$
\det_{T_X\otimes E^*}\big({\,}^T\Theta_{E,h_t}+
(1-t)\alpha\,\omega_0\otimes\Id_{E^*}\big)^{1/r}=
f_t\,\omega_0^n,\qquad f_t>0,
\leqno(2.9)
$$
and trace free Hermite-Einstein conditions
$$
\omega_t^{n-1}\wedge\Theta^\circ_{E,h_t}=g_t,\leqno(2.9^\circ)
$$
with smoothly varying families of functions $f_t\in C^\infty(X,\bR)$,
Hermitian metrics $\omega_t>0$ on~$X$ and sections
$g_t\in C^\infty(X,\Lambda^{n,n}_\bR T_X^*\otimes\Herm_{h_t}^\circ(E,E))$,
$t\in[0,1]$. Here, we start e.g.\ with the Yau-Uhlenbeck
solution $h_0=q_\varepsilon$ of (2.6) (so that $\det h_0=\det H_0$),
and take $\alpha>0$ so large that
${\,}^T\Theta_{E,h_0}+\alpha\,\omega_0\otimes\Id_{E^*}>0$ in the sense
of Nakano. If these conditions can be met for all $t\in[0,1]$
without any explosion of the solutions $h_t$, we infer from (2.9) that
$$
{\,}^T\Theta_{E,h_t}+
(1-t)\alpha\,\omega_0\otimes\Id_{E^*}>0
\quad\hbox{in the sense of Nakano} \leqno(2.9^+)
$$
for all $t\in[0,1]$. At time $t=1$, we will then get a Hermitian
metric $h_1$ on $E$ such that $\Theta_{E,h_1}$ is dual Nakano
positive. We still have the freedom of adjusting $f_t$, $\omega_t$
and $g_t$ in equations~$(2.9)$ and $(2.9^\circ)$. We have
a system of differential equations of order $2$, and any choice
of the right hand sides of the form
$$
\leqalignno{
&f_t(z)=F(t,z,h_t(z),D_zh_t(z),D^2_zh_t(z))>0,&(2.10)\cr
\noalign{\vskip6pt}
&g_t(z)=G(t,z,h_t(z),D_zh_t(z),D^2_zh_t(z))
\in C^\infty(X,\Lambda^{1,1}_\bR T^*_X\otimes\Herm^\circ(E,E))&(2.10^\circ)
\cr}
$$
is a priori acceptable for the sake of enforcing the positivity
condition $(2.9^+)$, although the presence of second order terms
$D^2_zh_t(z))$ might affect the principal symbol of the equations.
In equation $(2.9^\circ)$, the metrics
$\omega_t$ could possibly be taken to depend on~$t$, but unless some
commodity reason would appear in next stages of the analysis, it seems
simpler to set $\omega_t=\omega_0$ independent of~$t$. 
At this stage, we have the following\medskip

\noindent
{\bf 2.11. Theorem.} {\it Let $(E,H_0)$ be a smooth Hermitian
holomorphic vector bundle such that $E$ is ample and
$\omega_t=\omega_0=\Theta_{\det E,\det H_0}>0$.
Then the system of equations $(2.9,\;2.9^\circ)$ is a well determined
$($essentially non linear$)$ elliptic system of equations for all
choices of smooth right hand sides
$$
f_t=F(t,z,h_t,D_zh_t)>0,\qquad
g_t=G(t,z,h_t,D_zh_t,D^2_zh_t)\in\Herm^\circ(E,E),
$$
provided that the symbol $\eta_{h}$ of the linearized operator
$u\mapsto DG_{D^2h}(t,z,h,Dh,D^2h)\cdot D^2u$ has an
Hilbert-Schmidt norm
$\sup_{\xi\in T^*_X,|\xi|_{\omega_0}=1}\Vert\eta_{h}(\xi)\Vert_h\le
(r^2+1)^{-1/2}\,n^{-1}$ for any of the metrics $h=h_t$ involved. If a 
smooth solution $h_t$ exists on the whole time interval $[0,1]$,
then $E$ is dual Nakano positive.}\medskip\penalty-10000

\noindent
{\it Proof.} If we write a hermitian metric $h$ on $E$ under the form
\hbox{$h(v,w)=\langle\,\widetilde h(v),w\rangle_{H_0}$} with
$\widetilde h\in\Herm_{h_0}(E,E)$, we have
$h=H_0\smash{\widetilde h}$ in terms of
matrices. The curvature tensor is given by
the usual formula $\Theta_{E,h}=i\ol\partial(h^{-1}\partial h)=
i\ol\partial(\widetilde h^{-1}\partial_{H_0}\widetilde h)$, where
$\partial_{H_0}s=H_0^{-1}\partial(H_0s)$ is the $(1,0)$-component of
the Chern connection associated with~$H_0$ on~$E$. For simplicity
of notation, we~put
$$
M:=\Herm(E),\qquad M_h=\Herm_h(E,E),\qquad
M_h^\circ=\Herm_h^\circ(E,E).
$$
The system of equations $(2.9,2.9^\circ)$ is associated with the non linear
differential operator
$$
P:C^\infty(X,M)\to C^\infty(X,\bR\oplus M_h^\circ),\qquad
h\mapsto P(h)
$$
defined by
$$
P(h)=\omega_0^{-n}\Big(\!\det_{T_X\otimes E^*}\big({\,}^T\Theta_{E,h}+
(1-t)\alpha\,\omega_0\otimes\Id_{E^*}\big)^{1/r},\,
\omega_0^{n-1}\wedge\Theta_{E^\circ,h}-G(t,z,h,Dh, D^2h)\Big).
$$
It is by definition elliptic at $h$ 
if its linearization $u\mapsto (dP)_h(u)$ is an elliptic linear operator,
a crucial fact being that $M$ and \hbox{$\bR\oplus M_h^\circ$} have the same
rank $r^2$ over the field~$\bR$. Our goal is to compute the symbol
$\sigma_{dP}\in C^\infty(X,S^2T^\bR_X\otimes\Hom(M,\bR\oplus M_h^\circ))$
of $dP$, and to check that \hbox{$u\mapsto \sigma_{dP}(\xi)\cdot u$} 
is invertible for every non zero vector $\xi\in T^*_X$.
We  pick an infinitesimal varia\-tion
$\delta h$ of $h$ in $C^\infty(X,M)$, and represent it as
$\delta h=\langle u\,\bu,\bu\rangle_h$ with $u\in M_h=\Herm_h(E,E)$.
In terms of matrices, we have $\delta h=hu$,
i.e.\ $u=(u_{\lambda\mu})=h^{-1}\delta h$ is
the ``logarithmic variation of~$h$''. In~this setting, we evaluate
$(dP)_h(u)$ in orthonormal coordinates $(z_j)_{1\le j\le n}$ on $X$ rela\-tively
to~$\omega_0$. We have $h+\delta h=h(\Id+u)$ and
$(h+\delta h)^{-1}=(\Id-u)h^{-1}$ modulo $O(u^2)$, thus
$$
\leqalignno{
\qquad d\Theta_{E,h}(u)
&=i\ol\partial(h^{-1}\partial(hu))-i\ol\partial(uh^{-1}\partial h)
=i\ol\partial\partial u+i\ol\partial(h^{-1}\partial h\,u)
-i\ol\partial(uh^{-1}\partial h)\cr
&=i\ol\partial\partial_{h^*\otimes h} u=-i\partial_{h^*\otimes h}\ol\partial u,
&(2.12)}
$$
where $\partial_{h^*\otimes h}$ denotes here the $(1,0)$-component of the Chern
connection on the holomorphic vector bundle $\Hom(E,E)=E^*\otimes E$ 
induced by the metric $h^*\otimes h$.
As a consequence, the order $2$ term of the linearized operator is just
$$
d\Theta_{E,h}(u)^{[2]}=-i\partial\ol\partial u,
$$
and the logarithmic differential of the first scalar component
$P_\bR(h)$ of $P(h)$ has order $2$ terms given by
$$
P_{\bR}(h)^{-1}\,dP_{\bR,h}(u)^{[2]}={1\over r}
\tr(-\theta^{-1}\cdot {\,}^T i\partial\overline\partial u)=
-{1\over r}\,(\det\theta)^{-1}\sum_{j,k,\lambda,\mu}\widetilde\theta_{jk\lambda\mu}\,
{\partial^2 u_{\lambda\mu}\over\partial z_j\partial\ol z_k},\leqno(2.13)
$$
where $\theta$ is the $(n\times r)$-matrix of
$\theta=\theta(t,h)={\,}^T\Theta_{E,h}+(1-t)\alpha\omega_0\otimes\Id_{E^*}>0$,
$\;\widetilde\theta$ its co-adjoint and $\theta^{-1}=(\det\theta)^{-1}{\,}^T
\widetilde\theta$, so that $P_\bR(h)=\omega_0^{-n}\,(\det\theta)^{1/r}$.
We~also have to compute the order~$2$ terms in the differential of
the second component
$$
h\mapsto P^\circ(h)=\omega_0^{-n}\big(
\omega_0^{n-1}\wedge\Theta_{E,h}^\circ-G(t,z,h,Dh, D^2h)\big).
$$
Let us set $u={1\over r}\tr u\otimes\Id_E+u^\circ$,
$u^\circ\in M^\circ$, and $\tr u=\sum_\lambda u_{\lambda\lambda}\in\bR$.
Putting $\tau={1\over r}\tr u$, this actually gives an isomorphism
$M_h\to\bR\oplus M^\circ_h$, $u\mapsto(\tau,u^\circ)$.
Since $u^\circ$ is the logarithmic variation of $h^\circ=h(\det h)^{-1/r}$,
we get
$$
(dP^\circ)_h(u)^ {[2]}=\omega_0^{-n}\big(
-\omega_0^{n-1}\wedge i\partial\ol\partial u^\circ
-DG_{D^2h}\cdot D^2u\big).\leqno(2.14)
$$
If we fix a Hermitian metric $h$ and take a non zero 
cotangent vector $0\neq\xi\in T_X^*$, the symbol $\sigma_{dP}$ is given by
an expression of the form
$$
\sigma_{(dP)_h}(\xi)\cdot u= -\Bigg(
{(\det\theta)^{-1+1/r}\over r\,\omega_0^n}\sum_{j,k,\lambda,\mu} 
\tilde\theta_{jk\lambda\
mu}\,\xi_j\overline\xi_k\,u_{\lambda\mu}\,,\,
{1\over n}\,|\xi|^2u^\circ+\widetilde\sigma_G(\xi)\cdot u
\Bigg)\leqno(2.15)\phantom{\raise5pt\hbox{$\Big|$}}
$$
where $\widetilde\sigma_G$ is the principal
symbol of the operator $DG_{D^2h}\cdot D^2$.
If $g_t=G(t,z,h_t,Dh_t)$ is independent of $D^2h_t$,
the latter symbol $\widetilde\sigma_G$ is equal to~$0$ and
it is easy to see from (2.13) that $u\mapsto \sigma_{(dP)_h}(\xi)\cdot u$ is
an isomorphism in $\Hom(M_h,\bR\oplus M^\circ_h)$.
In fact, the first summation yields
$$
\sum_{j,k,\lambda,\mu} 
\tilde\theta_{jk\lambda\mu}\,\xi_j\overline\xi_k\,u_{\lambda\mu}=
\sum_{j,k,\lambda,\mu} 
\tilde\theta_{jk\lambda\mu}\,\xi_j\overline\xi_k\,u^\circ_{\lambda\mu}
+{1\over r}\sum_{j,k,\lambda} 
\tilde\theta_{jk\lambda\lambda}\,\xi_j\overline\xi_k\,\tr u.
$$
By an easy calculation, we get an inverse
operator $\bR\oplus M_h^\circ\to M_h$, $(\tau,v)\mapsto u$ where
$$
-r\,\omega_0^n\,(\det\theta)^{1-1/r}\,\tau=
\sum_{j,k,\lambda,\mu} 
\tilde\theta_{jk\lambda\mu}\,\xi_j\overline\xi_k\,u^\circ_{\lambda\mu}
+{1\over r}\sum_{j,k,\lambda} 
\tilde\theta_{jk\lambda\lambda}\,\xi_j\overline\xi_k\,\tr u,\quad
-v={1\over n}\,|\xi|^2u^\circ,
$$
hence $u^\circ=-{n\over|\xi|^2}v$ and
$$
\sigma_{(dP)_h}(\xi)^{-1}\cdot(\tau,v)=
{{n\over|\xi|^2}\sum_{j,k,\lambda,\mu}\tilde\theta_{jk\lambda\mu}\,
\xi_j\overline\xi_k\,v_{\lambda\mu}-r\,\omega_0^n\,(\det\theta)^{1-1/r}\,\tau
\over \sum_{j,k,\lambda}\tilde\theta_{jk\lambda\lambda}\xi_j\ol\xi_k}\;\Id_E
-{n\over|\xi|^2}\,v.
$$
Let us take the Hilbert-Schmidt norms 
$|u|^2=\smash{\sum_{\lambda,\mu}|u_{\lambda\mu}|^2}$
on $M_h=\Herm_h(E,E)$, and\break $c|\tau|^2+|v|^2$ on $\bR\oplus M_h^\circ$ ($h$
being the reference metric, and $C>0$ a constant). By homo\-geneity, we can
also assume \hbox{$|\xi|=|\xi|_{\omega_0}=1$}. Since
$(\sum_{j,k}\tilde\theta_{jk\lambda\mu}\,\xi_j\overline\xi_k)_{1\le \lambda,\mu\le r}$
is a positive Hermitian matrix by the Nakano positivity pro\-perty,
its trace is a strict upper bound for the largest eigenvalue, and we get
$$
\Bigg|\sum_{j,k,\lambda} 
\tilde\theta_{jk\lambda\mu}\,\xi_j\overline\xi_k\,v_{\lambda\mu}\Bigg|^2
\le (1-\delta)\Bigg(
\sum_{j,k,\lambda}\tilde\theta_{jk\lambda\lambda}\xi_j\ol\xi_k\Bigg)^{\!2}
\;\sum_\lambda|v_{\lambda\mu}|^2.
$$
The Cauchy-Schwarz inequality yields
$$
\Bigg|\sum_{j,k,\lambda,\mu} 
\tilde\theta_{jk\lambda\mu}\,\xi_j\overline\xi_k\,v_{\lambda\mu}\Bigg|^2
\le r(1-\delta)\,
\Bigg(\sum_{j,k,\lambda}\tilde\theta_{jk\lambda\lambda}\xi_j\ol\xi_k
\Bigg)^{\!2}\;\sum_{\lambda,\mu}|v_{\lambda\mu}|^2.
$$
For $|\xi|=1$, as $\Id_E\perp M^\circ$ and $|\Id_E|^2=r$, this implies
$$
\eqalign{
\big|\sigma_{(dP)_h}(\xi)^{-1}\cdot(\tau,v)\big|^2
&\le\Bigg(nr^{1/2}(1-\delta)^{1/2}|v|+{r\,\omega_0^n\,(\det\theta)^{1-1/r}
\over \sum_{j,k,\lambda}\tilde\theta_{jk\lambda\lambda}\xi_j\ol\xi_k}\,|\tau|
\Bigg)^{\!2}r+n^2\,|v|^2\cr
\noalign{\vskip5pt}
&< (n^2r^2+n^2)(C|\tau|^2+|v|^2)\cr}
$$
for $C$ large enough. By a standard pertubation argument, (2.13) remains
bijective if $|\widetilde\sigma_G(\xi)|_h$ is less than the inverse
of the norm of $\sigma_{(dP)_h}(\xi)^{-1}$, i.e.\ $(r^2+1)^{-1/2}\,n^{-1}$.
Similarly, one could also allow the scalar right
hand side $F$ to have a ``small dependence''
on~$D^2h_t$, but this seems less useful.\hfill\qed\medskip

Our next concern is to ensure that the existence of solutions holds on an
open interval of time $[0,t_0[$ (and hopefully on the whole interval $[0,1]$).
In~the case of a rank one metric $h=e^{-\varphi}$, it is well-known
that the K\"ahler-Einstein equation
$(\omega_0+i\partial\ol\partial\varphi_t)^n=e^{tf+\lambda\varphi_t}\omega_0^n$
more easily results in getting openness and closedness of solutions
when applying the continuity method for $\lambda>0$, as the linearized
operator $\psi\mapsto \Delta_{\omega_{\varphi_t}}\psi-\lambda \psi$
is~always invertible.
One way to generalize the K\"ahler-Einstein condition to the case of
higher ranks $r\ge 1$ is to take
$$
f_t(z)=(\det H_0(z)/\det h_t(z))^\lambda\,a_0(z),\qquad\lambda\ge 0,
\leqno(2.16)
$$
where $a_0(z)=\omega_0^{-n}\det({}^T\Theta_{E,h_0}
+\alpha\omega_0\otimes\Id_{E^*})^{1/r}>0$ is chosen so that the equation is
satisfied by $h_0$ at $t=0$ (the choice $\lambda>0$ has the interest that
$f_t$ gets automatically rescaled by multiplying $h_t$ by a constant,
thus ensuring strict invertibility). For the trace free part,
what is needed is to introduce a friction term $g_t$ that helps
again in getting invertibility of the linearized operator, and could
possibly avoid the explosion of solutions when $t$ grows to~$1$.
A choice compatible with the Yau-Uhlenbeck solution (2.8) at $t=0$ is to take
$$
g_t=-\varepsilon\,(\det H_0(z)/\det h_t(z))^\mu\,
\omega_0^n\otimes\log\widetilde h_t^\circ,
\qquad\varepsilon>0,~\mu\in\bR,
\leqno(2.16^\circ)
$$
if one remembers that $\det h_0=\det H_0$.
These right hand sides do not depend on higher derivatives of $h_t$,
so Theorem~2.11 ensures the ellipticity of the differential system.
Moreover:\medskip

\noindent
{\bf 2.17. Theorem.} {\it For $\varepsilon\ge\varepsilon_0(h_t)$ and
$\lambda\ge\lambda_0(h_t)(1+\mu^2)$ with $\varepsilon_0(h_t)$ and
$\lambda_0(h_t)$ large enough, the elliptic differential
system defined by
$(2.9,~2.9^\circ)$ and $(2.16,~2.16^\circ)$, namely
$$
\eqalign{
&\omega_0^{-n}\,\det_{T_X\otimes E^*}\big({\,}^T\Theta_{E,h_t}+
(1-t)\alpha\,\omega_0\otimes\Id_{E^*}\big)^{1/r}=
\bigg({\det H_0(z)\over\det h_t(z)}\bigg)^\lambda\,a_0(z)\cr
&\omega_0^{-n}\,\big(\omega_0^{n-1}\wedge\Theta^\circ_{E,h_t}\big)=
-\varepsilon\,\bigg({\det H_0(z)\over\det h_t(z)}\bigg)^\mu\,
\log\widetilde h^\circ_t,\cr}
$$
possesses an invertible
elliptic linearization. As a consequence, for such values of
$\varepsilon$ and $\lambda$, there exists an
open interval $[0,t_0)\subset [0,1]$ on which the solution
$h_t$ exists.}\medskip

\noindent
{\it Proof.} We replace the operator
$P:C^\infty(X,M)\to C^\infty(X,\bR\oplus M_h^\circ)$ used in the proof
of Theorem~2.9 by $\smash{\widetilde P=(\widetilde P_\bR,\widetilde P^\circ)}$
defined by
$$
\eqalign{
\widetilde P_\bR(h)&=\omega_0^{-n}(\det h(z)/\det H_0(z))^\lambda\,
\det_{T_X\otimes E^*}\big({\,}^T\Theta_{E,h}+
(1-t)\alpha\,\omega_0\otimes\Id_{E^*}\big)^{1/r},\cr
\noalign{\vskip5pt}
\widetilde P^\circ(h)&=
\omega_0^{-n}\big(\omega_0^{n-1}\wedge\Theta^\circ_{E,h}\big)+\varepsilon\,
(\det h(z)/\det H_0(z))^{-\mu}\,\log\widetilde h^\circ.\cr}
$$
Here, we have to care about the linearized operator $dP$ itself,
and not just with its principal symbol. We let again
$u=h^{-1}\delta h\in \Herm_h(E,E)$ and use formula (2.12)
for $d\Theta_{E,h}(u)$. This implies
$$
\widetilde P_{\bR}(h)^{-1}\,d\widetilde P_{\bR,h}(u)=\lambda\,\tr u
-{1\over r}\,\tr_{T_X\otimes E^*}\Big(\theta^{-1}\cdot{\,}^T\big(
i\partial_{h^*\otimes h}\ol\partial u\big)\Big).
$$
We need the fact that
$h^\circ=h\cdot(\det h)^{-1/r}$ possesses, when viewed as a
Hermitian endomorphism, a~logarithmic variation
$$
(\widetilde h^\circ)^{-1}\delta \widetilde h^\circ= u^\circ
= u -{1\over r}\tr u\cdot\Id_E.
$$
By the classical formula expressing the differential of
the logarithm of a matrix, we have
$$
d\log g(\delta g)=
\int_0^1\big((1-t)\Id+tg\big)^{-1}\delta g\,\big((1-t)\Id+tg\big)^{-1}\,dt,
$$
which implies
$$
d\log\widetilde h^\circ(u)=
\int_0^1\big((1-t)\Id+t\,\widetilde h^\circ\big)^{-1}
\,\widetilde h^\circ u^\circ\,\big((1-t)\Id+t\,\widetilde h^\circ\big)^{-1}\,dt.
$$
In the end, we obtain
$$
\eqalign{
&(d\widetilde P^\circ)_h(u)=
-\omega_0^{-n}\Big(\omega_0^{n-1}\wedge
i\partial_{h^*\otimes h}\ol\partial u^\circ\Big)+{}\cr
&~\varepsilon\,\bigg({\det h(z)\over \det H_0(z)}\bigg)^{-\mu}
\Bigg(\int_0^1\big((1-t)\Id+t\,\widetilde h^\circ\big)^{-1}
\,\widetilde h^\circ u^\circ\,
\big((1-t)\Id+t\,\widetilde h^\circ\big)^{-1}\,dt
-\mu\,\tr u\,\log\widetilde h^\circ\Bigg).\cr}
$$
In order to check the invertibility, we use the norm
$|\tau|^2+C|v|^2$ on $\bR\oplus M_h^\circ$ and compute the $L^2$ inner product
$\Langle(d\widetilde P)_h(u), (\tau,u^\circ)\Rangle$ over~$X$, where
$\tau={1\over r}\tr u$. The ellipticity
of operators $-i\partial_H\ol\partial$ implies that it has a discrete
sequence of eigenvalues converging to~$+\infty$, and that we get G{\aa}rding
type inequalities of the form  $\Langle -i\partial_H\ol\partial v,v\Rangle_H
\ge c_1\Vert\nabla v\Vert_H^2-c_2\Vert v\Vert_H^2$ where $c_1,c_2>0$ 
depend on~$H$. We apply such inequalities to
$v=\tau$, $H=1$, and $v=u^\circ$, $H=h^*\otimes h$, replacing
$u$ with $u=\tau\,\Id+u^\circ$. From~this, we infer
$$
\eqalign{
\Langle(d\widetilde P)_h(u), (\tau,u^\circ)\Rangle
&\ge c_1\Vert d\tau\Vert^2-c_2\Vert\tau\Vert^2 +\lambda r\,\Vert\tau\Vert^2
-{1\over r}\Langle\tr_{T_X\otimes E^*}\Big(\theta^{-1}\kern-1pt\cdot{}^T\big(
i\partial_{h^*\otimes h}\ol\partial u^\circ\big)\Big),\tau\Rangle\cr
&\qquad{}+C\Big(c_1^\circ\,\Vert\nabla u^\circ\Vert^2-c_2^\circ\Vert u^\circ\Vert^2
+c_3\varepsilon\,\Vert u^\circ\Vert^2-c_4\varepsilon\,|\mu|\,\Vert\tau\Vert\,
\Vert u^\circ\Vert\Big)\cr}
$$
where all constants $c_j$ may possibly depend on~$h$.
An integration by parts yields
$$
\eqalign{
\Big|{1\over r}\Langle\tr_{T_X\otimes E^*}\Big(\theta^{-1}\kern-1pt\cdot{}^T\big(
i\partial_{h^*\otimes h}\ol\partial u^\circ\big)\Big),\tau\Rangle\Big|
&\le c_5\Vert\nabla u^\circ\Vert\,(\Vert d\tau\Vert+\Vert\tau\Vert)\cr
&\le {1\over 2}c_1\,\big(\Vert d\tau\Vert^2+\Vert\tau\Vert^2\big)
+c_6\Vert\nabla u^\circ\Vert^2\cr}
$$
and we have
$$
c_4\varepsilon\,|\mu|\,\Vert\tau\Vert\,\Vert u^\circ\Vert\le
{1\over 2}c_3\varepsilon\,\Vert u^\circ\Vert^2+
c_7\varepsilon\mu^2\,\Vert\tau\Vert^2.
$$
If we choose $\varepsilon\ge 2c_2^\circ/c_3+1$, $C\ge c_6/c_1^\circ+1$ and
$\lambda r\ge c_2+{1\over 2}c_1+Cc_7\varepsilon\mu^2+1$, we finally get
$$
\Langle(d\widetilde P)_h(u), (\tau,u^\circ)\Rangle\ge
{1\over 2}c_1\Vert d\tau\Vert^2+\Vert\tau\Vert^2
+c_1^\circ\Vert\nabla u^\circ\Vert^2
+{1\over 2}Cc_3\varepsilon\Vert u^\circ\Vert^2
$$
and conclude that $(d\widetilde P)_h$ is an invertible elliptic operator.
The openness property at $t=0$ then follows from standard results on
elliptic PDE's.\hfill\qed
\medskip

\noindent
{\bf 2.18. Remarks.} (a) Theorem 2.17 is not very satisfactory since 
the constants
$\varepsilon_0(h_t)$ and $\lambda_0(h_t)$ depend on the solution $h_t$.
It would be important to know if one can get sufficiently uniform
estimates to make these constants independent of~$h_t$, thereby 
guaranteeing the long time existence of solutions. 
This might require
modifying somewhat the right hand side of our equations, especially
the trace free part, while taking a similar determinantal Monge-Amp\`ere
equation that still enforces the dual Nakano positivity of the curvature
tensor. The Yau iteration technique used in [Yau78] to get 0 order
estimates for Monge-Amp\`ere equations will probably have to be
adapted to this situation.

\noindent
(b) The non explosion of solutions when $t\to 1$  does not come
for free, since this property cannot hold when $\det E$ is ample,
but $E$ is not. One possibility would be
to show that an explosion at time $t_0<1$ produces a ``destabilizing 
subsheaf'' $\cS$ contradicting the ampleness of $E/\cS$, similarly
to what was done in [UhY86] to contradict the stability hypothesis.
\medskip

\noindent
{\bf 2.19. Variants.} (a) The determinantal equation always yields a
K\"ahler metric
$$
\beta_t:=\tr_E\big(\Theta_{E,h_t}+(1-t)\alpha\,\omega_0\otimes\Id_E\big)
=\Theta_{\det E,\det h_t}+r(1-t)\alpha\,\omega_0>0.
$$
An interesting variant of the trace free equation is
$$
\omega_0^{-n}\,\big(\omega_t^{n-1}\wedge\Theta^\circ_{E,h_t}\big)=
-\varepsilon\,\bigg({\det H_0(z)\over\det h_t(z)}\bigg)^\mu\,
\log\widetilde h^\circ_t\leqno(*)
$$
with $\omega_t={1\over r\alpha+1}\beta_t$
(notice that $\beta_0=(r\alpha+1)\omega_0$). It is then important to
know whether the corresponding differential system is still elliptic
with an invertible linearization. According to equation $(*)$,
the $\Herm(E,E)^\circ$ part of the differential system depends on
the functional
$$
\widetilde P^\circ(h)=
\omega_0^{-n}\big(\omega_t^{n-1}\wedge\Theta^\circ_{E,h}\big)+\varepsilon\,
(\det h(z)/\det H_0(z))^{-\mu}\,\log\widetilde h^\circ,
$$
and, with respect to the functional used in Theorem 2.17,
the differential $d\widetilde P^\circ_h(u)$ acquires one additional term
coming from the variation of $\omega_t^{n-1}$. With the same notation
as in our previous calculations, we have 
$\Theta_{\det E,\det h_t}=-i\partial\ol\partial\log\det(h_t)$ and
$\delta(\beta_t)_h(u)=-i\partial\ol\partial \tr u$, hence
$$
\eqalign{
&(d\widetilde P^\circ)_h(u)=
-\omega_0^{-n}\Big(\omega_t^{n-1}\wedge
i\partial_{h^*\otimes h}\ol\partial u^\circ
+{\textstyle{n-1\over r\alpha+1}}\,
\omega_t^{n-2}\wedge i\partial\ol\partial \tr u
\wedge\Theta^\circ_{E,h}\Big)+{}\cr
&~\varepsilon\,\bigg({\det h(z)\over \det H_0(z)}\bigg)^{-\mu}
\Bigg(\int_0^1\big((1-t)\Id+t\,\widetilde h^\circ\big)^{-1}
\,\widetilde h^\circ u^\circ\,
\big((1-t)\Id+t\,\widetilde h^\circ\big)^{-1}\,dt
-\mu\,\tr u\,\log\widetilde h^\circ\Bigg).\cr}
$$
Putting again $\tau=\tr u$, this requires to estimate one
extra term appearing in the $L^2$ inner product
$\Langle(d\widetilde P)_h(u), (\tau,u^\circ)\Rangle$, namely
$$
\Langle
(\omega_0^n)^{-1}\omega_t^{n-2}\wedge i\partial\ol\partial \tau
\wedge\Theta^\circ_{E,h}\,,\,u^\circ\Rangle.
$$
We can apply the same integration by part argument as before
to conclude that $(d\widetilde P)_h$ is again invertible,
under a similar hypothesis $\lambda\geq \lambda_0(h_t)(1+\mu^2)$,
at least for $t$ small. A~very recent note posted by V.P.\ Pingali 
[Pin21] shows that when $E$ is $\omega_0$-stable and $h_0$ is
taken to be the Hermite-Einstein metric, the trace free part of
the differential system used in Theorem 2.17 has a solution of
the form $h_t=h_0e^{-\psi_t}$, thus always ``conformal'' to $h_0$.
There are cases where the dual Nakano positivity of $h_0$ is doubtful.
As a consequence, even in that favorable case, it is 
unclear whether a long time existence result can hold for the
total system, unless stronger 
restrictions on the Chern classes are made. Equation $(*)$ does not
seem to entail such contraints, and may thus be better suited to the
investigated problem.

\noindent(b) In a first step towards solving (2.6), [UhY86] consider
equations that have even stronger friction terms, taking the right
hand side to be of the form
$$
\omega_0^{n-1}\wedge\Theta_{E,h}=\omega_0^n\otimes
\big(-\varepsilon\,\log\widetilde h+\sigma\,\widetilde h^{-1/2}\,
\Gamma_0\,\widetilde h^{1/2}-\Gamma_0),\quad\sigma>0,
$$
and letting $\sigma\to 0$ at the end of the analysis. Here we can do
just the same, for instance by adding a term equal to a multiple of
$(\widetilde h^\circ_t)^{-1/2}\,\Gamma_t\,(\widetilde h^\circ_t)^{1/2}-\Gamma_t$
in the trace free equation, as such terms are precisely trace free for any
$\Gamma_t\in C^\infty(X,\Hom(E,E))$.
\bigskip

\noindent
{\bigbf 3. A concept of Monge-Amp\`ere volume for vector bundles}

If $E\to X$ is an ample vector bundle of rank $r$, the associated line
bundle
$$
\cO_{\bP(E)}(1)\to Y=\bP(E)
$$
is ample, and one can consider its volume
$c_1(\cO_{\bP(E)}(1))^{n+r-1}$. It is well known that this number (which is 
an integer) coincides with the Segre number $\int_X (-1)^ns_n(E)$, where
$(-1)^ns_n(E)$ is the $n$-th Segre class of~$E$. Let us assume further
that $E$ is dual Nakano positive (if the solution of the Hermitian-Yangs-Mills
differential system of \S2 is unobstructed, this would follow from the 
ampleness of $E$). One
can then introduce the following more involved concept of volume,
which we will call the Monge-Amp\`ere volume of $E\,$:
$$
\MAVol(E)=\sup_h\int_X\det_{T_X\otimes E^*}\big(
(2\pi)^{-1}{\,}^T\Theta_{E,h}\big)^{1/r},
\leqno(3.1)
$$
where the supremum is taken over all smooth metrics $h$ on $E$ such that
${}^T\Theta_{E,h}$ is Nakano positive. This supremum is always finite, and 
in fact we have\medskip

\noindent
{\bf 3.2. Proposition.} {\it For any dual Nakano positive vector bundle $E$, 
one has
$$\MAVol(E)\leq r^{-n}c_1(E)^n.$$}%%
\noindent{\it Proof.} We take $h$ to be a hermitian metric on $E$ such that
${}^T\Theta_{E,h}$ is Nakano positive, and consider the K\"ahler metric
$$
\omega=(2\pi)^{-1}\Theta_{\det E,\det h}=
(2\pi)^{-1}\tr_{E^*}{}^T\Theta_{E,h}\in c_1(E).
$$
If $(\lambda_j)_{1\le j\le nr}$ are the eigenvalues of the associated
hermitian form $(2\pi)^{-1}{}^T\widetilde\Theta_{E,h}$ with respect to
$\omega\otimes h$, we have
$$
\det_{T_X\otimes E^*}\big(
(2\pi)^{-1}{\,}^T\Theta_{E,h}\big)^{1/r}
=\bigg(\prod_j\lambda_j\bigg)^{1/r}\omega^n
$$
and $\big(\prod_j\lambda_j\big)^{1/nr}\le
{1\over nr}\sum_j\lambda_j$ by the inequality between geometric and
arithmetic means. Since
$$
\sum_j\lambda_j=\tr_\omega\big(\tr_{E^*}\big((2\pi)^{-1}{\,}^T\Theta_{E,h}
\big)\big)=\tr_\omega\omega=n,
$$
we conclude
$$
\int_X \det_{T_X\otimes E^*}\big(
(2\pi)^{-1}{\,}^T\Theta_{E,h}\big)^{1/r}
\le \int_X \bigg({1\over nr}\sum_j\lambda_j\bigg)^n\omega^n
=r^{-n}\int_X\omega^n=r^{-n}c_1(E)^n.
$$
The proposition follows.\hfill\qed
\medskip

\noindent
{\bf 3.3. Remarks.} (a) In case $E=\bigoplus_{1\le j\le r}E_j$ and
$h=\bigoplus_{1\le j\le r}h_j$ are split, with all metrics $h_j$ normalized
to have proportional volume forms
$((2\pi)^{-1}\Theta_{E_j,h_j})^n=\beta_j\omega^n$ with suitable constants
$\beta_j>0$, we get $\beta_j=c_1(E_j)^n/c_1(E)^n$, and the inequality
reads
$$
\bigg(\prod_{1\le j \le r}c_1(E_j)^n\bigg)^{1/r}\le r^{-n}c_1(E)^n.
$$
It is an equality when $E_1=\cdots=E_r$, thus Proposition 3.2 is
optimal as far as the constant $r^{-n}$ is concerned.
For $E=\bigoplus_{1\le j\le r}E_j$ split with distinct ample factors,
it seems natural to conjecture
that
$$
\MAVol(E)=\bigg(\prod_{1\le j \le r}c_1(E_j)^n\bigg)^{1/r},
$$
i.e.\ that the supremum is reached for split metrics $h=\bigoplus h_j$.
In case $E$ is a non split extension $0\to A\to E\to A\to 0$ with
$A$ an ample line bundle -- this is possible if $H^1(X,\cO_X)\ne 0$,
e.g.\ on an abelian variety -- we strongly suspect that
$\MAVol(E)=c_1(A)^n$ but that the supremum  is not reached by any
smooth metric, as we have $E$ semi-stable  but not polystable.
\medskip

\noindent
(b) It would be interesting to characterize the ``extremal
metrics'' $h$ achieving the supremum in (3.1) when they exist. The
calculations made in \S2 show that they satisfy some
Euler-Lagrange equation
$$
\int_X
(\det \theta)^{1/r}\cdot
\tr_{T_X\otimes E^*}\Big(\theta^{-1}\cdot{\,}^T\big(
i\partial_{h^*\otimes h}\ol\partial u\big)\Big)=0\qquad
\forall u\in C^\infty(X,\Herm(E)),
$$
where $\theta$ is the $(n\times r)$-matrix representing
${\,}^T\Theta_{E,h}$. After performing two integration by parts
that free $u$ from any differentiation, we get a fourth order non linear
differential system that $h$ has to satisfy. Remark 3.3~(a) leads us 
to suspect that this system is not always solvable, but the
addition of adequate lower order ``friction terms'' might make it 
universally solvable. This could possibly yield a better alternative to the
more naive order 2 differential system we proposed in \S2 to study the
Griffiths conjecture.
\medskip
\noindent
(c) When $r>1$, one may wonder what is the infimum
$$
\inf_h\int_X\det_{T_X\otimes E^*}\big(
(2\pi)^{-1}{\,}^T\Theta_{E,h}\big)^{1/r}.
$$
In the split case $(E,h)=\bigoplus(E_j,h_j)$, we can normalize
$\Theta_{E_j,h_j}$ to satisfy $\Theta_{E_j,h_j}^n=f_j\omega^n$ with
$\int_X f_j\omega^n=c_1(E_j)^n$, $f_j>0$. Then
$$
\int_X\det_{T_X\otimes E^*}\big((2\pi)^{-1}{\,}^T\Theta_{E,h}\big)^{1/r}=
\int_X(f_1\cdots f_r)^{1/r}\,\omega^n
$$
and this integral becomes arbitray small if we take the
$f_j$'s to be large on disjoint open sets, and very small elsewhere.
This example leads us to suspect that one always has
$$
\inf_h\int_X\det_{T_X\otimes E^*}\big(
(2\pi)^{-1}{\,}^T\Theta_{E,h}\big)^{1/r}=0
$$
for $r>1$. The ``friction terms'' used in our differential systems
should be chosen so as to prevent any such shrinking of the volume.
\vskip30pt

\centerline{\bigbf References}
\vskip8pt

\Bibitem[Ber09]&Berndtsson B.:& Curvature of vector bundles associated to
holomorphic fibrations,& Annals of Math.\ {\bf 169} (2009), 531–-560&

\Bibitem [CaF90]&Campana F., Flenner, H.:& A characterization of ample
vector bundles on a curve,& Math.\ Ann.\ {\bf 287} (1990),  571--575&

\Bibitem [DeS79]&Demailly, J.-P., Skoda, H.:& Relations entre les
notions de positivit\'e de P.A. Griffiths et de S. Nakano,&
S\'eminaire P.~Lelong-H.~Skoda (Analyse), ann\'ee 1978/79, Lecture
notes in Math., no 822, Springer-Verlag, Berlin (1980), 304--309&

\Bibitem [Don85]&Donaldson, S.:& Anti-self-dual Yang-Mills connections over
complex algebraic surfaces and stable vector bundles,& Proc.\ London Math.\
Soc. {\bf 50} (1985), 1--26&

\Bibitem [Gri69]&Griffiths, P.A:& Hermitian differential geometry, Chern classes and positive vector bundles,& Global Analysis, papers in honor of K. Kodaira, Princeton Univ. Press, Princeton (1969), 181--251&

\Bibitem [Kod54]&Kodaira, K.:& On K\"ahler varieties of restricted type,&
Ann.\ of Math.\ {\bf 60} (1954) 28--48&

\Bibitem [MoT07]&Mourougane, C., Takayama, S.:& Hodge metrics and positivity of direct images,& J.~reine angew.\ Math.\ {\bf 606} (2007), 167--179&

\Bibitem [Nak55]&Nakano, S.:& On complex analytic vector bundles,& J.~Math.\ Soc.\ Japan {\bf 7} (1955) 1--12&

\Bibitem [NaS65]&Narasimhan, M. S., Seshadri, C. S.:& Stable and unitary
vector bundles on a compact Riemann surface,& Ann.\ of Math.\ {\bf 82} (1965),
540--567&

\Bibitem [Nau18]&Naumann, P.:& An approach to Griffiths conjecture,&
arXiv:1710.10034, math.AG&

\Bibitem [Pin20]&Pingali, V.P.:& A vector bundle version of the Monge-Amp\`ere
equation,& Adv.\ in Math.\ {\bf 360} (2020), 40 pages,
https://doi.org/10.1016/j.aim.2019.106921&

\Bibitem [Pin21]&Pingali, V.P.:& A note on Demailly's approach towards a
conjecture of Griffiths,& manuscript, Feb.~2021, to appear in
C.~R.\ Math.\ Acad.\ Sci.\ Paris&

\Bibitem [UhY86]&Uhlenbeck, K., Yau, S.T.:& On the existence of 
Hermitian-Yang-Mills connections in stable vector bundles,& Comm.\
Pure and Appl.\ Math.\ {\bf 39} (1986) 258--293&

\Bibitem [Ume73]&Umemura, H.:& Some results in the theory of vector bundles,&
Nagoya Math.\ J.\ {\bf 52} (1973), 97--128&

\Bibitem [Yau78]&Yau, S.T.:& On the Ricci curvature of a complex K\"ahler
manifold and the complex Monge-Amp\`ere equation I,& Comm.\ Pure and Appl.\
Math.\ {\bf 31} (1978) 339--411&
\vskip15pt

\noindent
Jean-Pierre Demailly\hfil\break
Universit\'e Grenoble Alpes, Institut Fourier (Math\'ematiques)\hfil\break
UMR 5582 du C.N.R.S., 100 rue des Maths, 38610 Gi\`eres, France\hfil\break
{\it e-mail:}\/ jean-pierre.demailly@univ-grenoble-alpes.fr

\bye